\documentclass[12pt]{elsart}
\usepackage[latin1]{inputenc}
\usepackage{amsmath, amsfonts, amssymb}
\usepackage{stmaryrd}
\usepackage{mathrsfs}
\usepackage{latexsym}

\newtheorem{theorem}{Theorem}[section]
\newtheorem{lemma}[theorem]{Lemma}
\newtheorem{proposition}[theorem]{Proposition}
\newtheorem{corollary}[theorem]{Corollary}

\newtheorem{definition}[theorem]{Definition}
\theoremstyle{definition}
\newtheorem{remark}[theorem]{Remark}

\newcommand{\HS}[1]{\mathcal{H}(#1)}

\newcommand{\N}{\mathbb{N}}

\newcommand{\R}{\mathbb{R}}

\newcommand{\U}{\mathcal{U}}
\newcommand{\UU}{\mathbb{U}}

\newcommand{\bigO}{\ensuremath{\mathrm{{{O}}}}} 

\newcommand{\complarg}[1]{\ensuremath{\mathrm{{{#1}}}}} 
\newcommand{\complK}{\complarg{K}}
\newcommand{\complC}{\complarg{C}}
\newcommand{\complKm}{\complarg{Km}}
\newcommand{\complKM}{\complarg{KM}}

\newcommand{\zo}{\{0,1\}}
\newcommand{\cs}{2^\omega}
\newcommand{\uh}{{\upharpoonright}}
\renewcommand{\phi}{\varphi}
\newcommand{\str}{2^{<\omega}}
\newcommand{\fs}{\str}

\newcommand{\leqx}{\preceq}

\newcommand{\tuple}[1]{\langle #1 \rangle}

\journal{Journal of Computer and System Sciences}

\begin{document}

\begin{frontmatter}



\title{Solovay functions and their applications in algorithmic randomness}
\maketitle


\author[label1]{Laurent Bienvenu}
\author[label2]{Rod Downey}
\author[label3]{Andr\'e Nies}
\author[label4]{Wolfgang Merkle}
\address[label1]{Laboratoire CNRS J.-V. Poncelet, Moscow, Russia}
\address[label2]{School of Mathematics, Statistics and Operations Research, Victoria University of Wellington, New Zealand}
\address[label3]{Department of Computer Science,  University of Auckland, New Zealand}
\address[label4]{Institut f\"ur Informatik, Heidelberg University, Germany}


\begin{keyword}
Kolmogorov complexity, Algorithmic randomness, Solovay functions. 
\end{keyword}

\begin{abstract}
Classical versions of Kolmogorov complexity are incomputable. Nevertheless, in 1975 Solovay showed that there are computable functions~$f \geq \complK+O(1)$ such that for infinitely many strings $\sigma$, $f(\sigma)=\complK(\sigma)+O(1)$, where $\complK$ denotes prefix-free
Kolmogorov complexity (while $\complC$ denotes plain Kolmogorov complexity). Such an $f$ is now called a Solovay function. We prove that many classical results about $\complK$ can be obtained by replacing $\complK$ by a Solovay function. For example, the three following properties of a function $g$ all hold for the function K.

(i) The sum of the terms $\sum_n 2^{-g(n)}$ is a Martin-L\"of random real.\\
(ii) A sequence A is Martin-L\"of random if and only if $\complC(A \upharpoonright n) > n -g(n)-O(1)$.\\
(iii) A sequence A is K-trivial if and only if $\complK(A \upharpoonright n) < g(n) + O(1)$.\\

\noindent We show that when fixing any of these three properties, then among all computable functions exactly the Solovay functions possess this property. Furthermore, this characterization extends accordingly to the larger class of right-c.e. functions.
\end{abstract}

\end{frontmatter}


\section{Introduction and overview}\label{sec:introduction}
\subsection{Introduction}
A fundamental aspect of Kolmogorov complexity is its inherent noncomputability. That is,
standard complexities measuring compressibility, such as plain (Kolmogorov) complexity 
$\complC(\sigma)$ or prefix-free (Kolmogorov) complexity $\complK(\sigma)$, 
are  as complicated to calculate as the halting problem, and hence undecidable.
This nonalgorithmic aspect of Kolmogorov complexity 
is useful in the sense that it can enable new proofs
of various undecidability results such as G\"odel's Incompleteness Theorems (see~\cite[Chapter 2]{LiV2008} for the First Incompleteness Theorem, and~\cite{KritchmanR2010} for the Second Incompleteness Theorem) . On the other hand, this aspect also hinders 
the use of Kolmogorov complexity as a tool for measuring common 
information in real data. This  hindrance necessitates 
the use of replacements for the classical complexities by
computable text compressions (see for example, Cilibrasi and Vit\'{a}nyi \cite{CilibrasiV2005}).

Nevertheless,  it is remarkable how much of the classical theory of 
Kolmogorov complexities can be carried out by using 
\emph{good effective  upper bounds}. For example, for plain complexity we know that
the identity machine describing $\sigma$ by  $\sigma$ itself will give us 
$f(\sigma)=|\sigma|$ as a simple computable upper bound, and Kolmogorov's 
basic counting argument shows that this bound is achieved  
for each $n$ by some $\sigma$ of length $n$. 
Thus the function $f(\sigma)=|\sigma|$ is a good computable 
upper bound for $\complC(\sigma)$ in the sense that 
\begin{itemize}
\item[(i)] $\complC(\sigma)\le f(\sigma)+\bigO(1)$, and 
\item[(ii)] $\exists^\infty \sigma\,  f(\sigma)\le\complC(\sigma)+\bigO(1)$. 
\end{itemize}
This upper bound has a number of uses. One 
is Chaitin's result~\cite{Chaitin1976} that 
if $A$ is $C-$\emph{trivial}, which means that 
$\complC(A\upharpoonright n)\le \complC(n) +\bigO(1)$ for all $n$, then $A$ is computable.

The situation for prefix-free complexity is not so clear. Similar to the simple upper bound $|\sigma|$ for $\complC(\sigma)$, the function~$|\sigma|+\complK(|\sigma|)$ is an upper bound for~$\complK(\sigma)$ that is infinitely often tight up to an additive constant. However, this upper bound is not computable but only right-c.e., i.e., the binary relation ``$g(n) < k$'' is c.e. Nevertheless, by a theorem of Solovay there are such infinitely often tight upper bounds for~$\complK$ that are computable and such functions are now called Solovay functions.

\begin{definition} \label{def:solovayfunction}
A function $g$ 
is a \emph{Solovay function} if~$g$ is computable and it holds that 
\begin{itemize}
\item[(i)] $\complK(\sigma)\le g(\sigma)+\bigO(1)$, and 
\item[(ii)] $\exists^\infty \sigma\ g(\sigma) \le \complK(\sigma)+\bigO(1)$. 
\end{itemize}
A function~$g$ is a \emph{weak Solovay function} if~$g$ is right-c.e.\ and satisfies~(i) and~(ii).   
\end{definition}

\bigskip
In what follows we derive a number of fundamental results about Solovay functions and weak Solovay functions. Recall that 
the three following properties of a function~$g$ all hold for the function~$\complK$, by well-known results or, in case of the third one, by definition. 
\begin{itemize}
\item[(i)] The sum~$\sum_{n} 2^{-g(n)}$ is a Martin-L\"{o}f random real. 
\item[(ii)] A sequence~$A$ is Martin-L\"{o}f random if and only if~$\complC(A\upharpoonright n)\ge n - g(n) - \mathrm{O}(1)$.
\item[(iii)] A sequence~$A$ is K-trivial if and only if~$\complK(A\upharpoonright n)\le g(n) + \mathrm{O}(1)$.
\end{itemize}

Our main results are that when fixing any of these three properties, then the property is true for all weak Solovay functions but indeed is false for all right-c.e.\ functions that are not weak Solovay functions. That is, among all right-c.e.\ functions exactly the weak Solovay functions possess this property. This characterization extends trivially to the special case of Solovay functions in the sense that among all computable functions exactly the Solovay functions possess the property under consideration.  Note that in case of the first property the latter characterization was shown in a conference article by the first two authors of this paper~\cite{BienvenuD2009}, whereas the corresponding characterization of weak Solovay functions was subsequently demonstrated by H\"{o}lzl, Kr\"{a}ling and Merkle~\cite{HolzlKM2013}, see Section~\ref{sec:hoelzletal} for further details.

Quite aside from their intrinsic interest, Solovay functions have found many uses in the study of algorithmic randomness, many using the results proven in this paper some of which had been reported in earlier conference papers~\cite{BienvenuD2009,BienvenuMN2011}. Recall that $A$ is called $\complK$-trivial if and only if  there is a constant $b\in \mathbb N$ such that $\complK(A\uh n)\leq  \complK(n)+b$ for all $n$. Barmpalias and Sterkenburg~\cite{BarmpaliasS2011} used  Solovay functions to show that calculating the number of $\complK$-trivials with constant~$d$ is sharply $\Delta_3^0$. More recently, Solovay functions were used for the characterization of $\complK$-trivial points in computable metric spaces by Melnikov and Nies \cite{MelnikovN2013}. A research announcement by Bienvenu, Day, et al.~\cite{BienvenuDayetal13}  contains a  new,  golden-run-free proof that $\complK$-triviality  implies lowness for $\complK$. One ingredient  is our new,  golden-run-free proof, based on Solovay functions, that every $\complK$-trivial is Turing below a c.e.\ $\complK$-trivial  (Section~\ref{s:deg K-trivs}). For background on the golden run method see Nies~\cite[Section 5.4]{Nies2009}.

\subsection{Overview}\label{sec:overview}
The properties~(i), (ii), and~(iii) of~$\complK$ stated in the last section, which characterize weak Solovay functions, are discussed in Section~\ref{sec:basics}, \ref{sec:mlrandomness}, and~\ref{sec:ktriviality}, respectively. In Section~\ref{sec:existence} we review Solovay's construction of a Solovay function before we discuss in Section~\ref{sec:hoelzletal} the characterization of weak Solovay functions by property~(i), i.e., by the fact that the sum of the terms~$2^{-g(n)}$ is a Martin-L\"{o}f random real. 
    
In Section~\ref{sec:mlrandomness}, for a start we investigate into conditions formulated in terms of complexity gaps that are sufficient or necessary for being Martin-L\"{o}f random. We derive several negative results in this direction, for example, there is no function~$h$ that tends to infinity such that $\complK(A\upharpoonright  n)\ge n-h  (n)-\bigO(1)$ implies that~$A$ is Martin-L\"of random. These results will be used in subsequent proofs but have some interest in their own. In particular, they contrast a positive result by Csima and Montalb{\'a}n~\cite{CsimaM2005}, who gave a sufficient condition for K-triviality in terms of a  complexity gap. Section~\ref{sec:gacsmilleryu} is then about characterizations of Martin-L\"{o}f randomness via weak Solovay functions and vice versa. Recall that the G{\'a}cs-Miller-Yu Theorem \cite{Gacs1980,MillerY2008} gives a characterization of Martin-L\"of randomness in terms of plain Kolmogorov complexity $\complC$: a sequence~$A$ is Martin-L\"of random if and only if $\complC(A \uh n) \geq n - \complK(n) - \bigO(1)$;  moreover, there exists a computable upper bound~$f$ for~$\complK$ such that $A$ is Martin-L\"of random if and only if $\complC(A \uh n) \geq n - f(n)-\bigO(1)$.
We show that any right-c.e.\ function~$f$ characterizes Martin-L\"of randomness by the latter condition if and only if~$f$ is a weak Solovay function.
					
In Section~\ref{sec:ktriviality}, we explore the relationship between $\complK$-triviality and Solovay functions, showing that for any weak Solovay function $g$, a sequence~$A$ is  $\complK$-trivial if and only if $\complK(A\upharpoonright n)\le g(n)+\bigO(1)$. Again the characterization works both ways. Namely, if $f$ is right-c.e.\ and has the property that $\complK(A\upharpoonright n)\le f(n)+\bigO(1)$ is equivalent to~$A$ being $\complK$-trivial, then $f$ must be a weak Solovay function, and a similar characterization holds for Solovay functions. Nies \cite{Nies2005} used the golden run method to prove that every $\complK$-trivial is truth-table reducible to a c.e.\ $\complK$-trivial. Using the connection  between Solovay functions and $\complK$-trivials, we get a new and easy proof of the slightly weaker result  that every $\complK$-trivial is Turing reducible to a  c.e.\ $\complK$-trivial.
  
Finally, in Section~\ref{sec:hitting}, we look at the $c$-\emph{hitting set}
of a Solovay function~$g$, which is the set of numbers~$n$ such that  $g(n)\le \complK(n)+c$. 
Any such set, unless finite, is shown to be Turing complete, and sparse, in the sense that it is hyperimmune.
                                                     
\subsection{Notation and preliminaries}\label{sec:notation}
Here we gather some notation that will be used throughout the paper. String refers to a binary string, and~$\fs$ denotes the set of all strings. The length of a string~$x$ is denoted by~$|x|$.
Both, plain and prefix-free Kolmogorov complexity are defined on the set of finite strings, but as usual are also applied to other objects such as integers, rational numbers, pairs of strings, etc., as long as they can be effectively encoded by finite strings. Unless explicitly stated otherwise, sequence refers to an infinite binary sequence. A subset~$A$ of the natural numbers may be identified with its characteristic sequence~$A(0)A(1) \ldots$.  Occasionally, sequences are referred to as \emph{reals}, see the discussion on page~\pageref{disc-reals}. The set of all sequences is denoted by $\cs$. For~$i=0,1, \ldots$, the $i$-th bit of a sequence~$A$ is denoted by $A(i)$, while the prefix of $A$ of length~$i$ is written~$A \uh i$, that is,~$A=A(0)A(1)\ldots$, and~$A \uh i=A(0) \ldots A(i-1)$ for~$i>0$, while~$A \uh 0$ is equal to the empty string. For a string~$\sigma$, the \emph{cylinder}~$[\sigma]$ is the set of sequences~$A$ such that $\sigma$ is a prefix of~$A$. If $S$ is a \emph{set} of strings, we write~$[S]$ for the set of sequences having some prefix in~$S$, i.e., $S=\bigcup_{\sigma \in S} [\sigma]$. We denote by~$\lambda$ the uniform measure on $\cs$, which is the probability measure one gets when each bit of a sequence is chosen at random with probabilities $(1/2,1/2)$ and independently of all the other bits.\\

An \emph{order} is a nondecreasing, unbounded function from $\N$ to $\N$. For a given order~$h$, let $h^{-1}(k)$ be the largest integer~$n$ such that $h(n) \leq k$. Note that for any order~$h$, the function~$h^{-1}$ is itself an order, and is computable if and only if~$h$ is.   \\

A function~$f\colon D \rightarrow \mathbb{R}$ is \emph{left-c.e.} (also known as \emph{approximable or semi-computable from below}) if there is a computable function~$(x,t) \mapsto f_t(x)$ such that for all~$x \in D$, $f_t(x)$ is nondecreasing in~$t$ and converges to~$f(x)$. The value~$f_t(x)$ is called the approximation of~$f(x)$ at stage~$t$. The notion of a \emph{right-c.e.} function is defined accordingly, where now the approximations~$f_t(x)$ are required to be nonincreasing in~$t$. 
Note that for any left-c.e.\ or right-c.e.\ function one can assume that the approximations~$f_t(x)$ are given as a pair of natural numbers~$p$ and~$q$, which represent the dyadic rational~$p 2^{-q}$. Furthermore, in the important case of right-c.e.\ functions with values in the natural numbers such as plain and prefix-free Kolmogorov complexity, by rounding approximations down to the next natural number below, one can assume that the approximations~$f_0(x), f_1(x), \ldots$ are natural numbers, and that hence for each~$x$ all but finitely many of the approximations~$f_t(x)$ are equal to~$f(x)$.\\

We will assume that the reader is familiar with the rudiments of algorithmic randomness, the definitions of plain and prefix-free complexity, and the like. We only remind the reader of some of the most salient points. A \emph{bounded request} set, also known as Kraft-Chaitin or KC-set, is a computably enumerable set~$W$ of pairs $(\sigma,n)$ where the first coordinate is a string and the second an integer and such that $\sum_{(\sigma,n) \in W} 2^{-n}$ is finite. Enumerating a pair $(\sigma,n)$ into a request set is often said to incur a \emph{cost} of $2^{-n}$; the request set being bounded refers to the fact that the total cost is finite. Given such a  bounded request set~$W$, the \emph{Kraft-Chaitin theorem}, due to Levin~\cite{Levin1971} and  Chaitin~\cite{Chaitin1975}, asserts that  $\complK(\sigma) \leq n +\bigO(1)$  for all $(\sigma,n) \in W$. 

The classical application of the Kraft-Chaitin Theorem is the Levin-Schnorr characterization of Martin-L\"of randomness. 

\begin{theorem}[Levin-Schnorr \cite{Levin1973,Schnorr1973}] 
\label{levin-s}
A sequence~$X$ is Martin-L\"of random if and only if $\complK(X \upharpoonright n)\geq n-\bigO(1)$.
\end{theorem}

Another consequence of the Kraft-Chaitin Theorem is that every right-c.e.\ function of finite weight is an upper bound for Kolmogorov complexity.

\begin{lemma}[Levin~\cite{Levin1971}, Chaitin~\cite{Chaitin1975}] \label{lem:kub-sum}
Let~$f\colon\N \rightarrow \N$ be a right-c.e.\ function. Then the two following assertions are equivalent.
\begin{itemize} 
\item[(i)] $\complK \leq f+\bigO(1)$.
\item[(ii)] $\sum_{n} 2^{-f(n)}$ is finite. 
\end{itemize}  
\end{lemma}

\begin{pf*}{Proof} 
The implication (i)$\Rightarrow$(ii) is trivial as $\sum_n 2^{-\complK(n)} \leq 1$. For the implication~(ii)$\Rightarrow$(i), fix some computable function~$(n,t)\mapsto f_t(n)$ with values in the natural numbers such that for all~$n$, the sequence~$f_0(n), f_1(n), \ldots$ converges nonincreasingly to~$f(n)$, and let the set~$D_n$ be equal to~$\{f_0(n), \ldots, f(n)\}$, hence~${\sum_{\ell \in D_n}  2^{-\ell} \leq 2^{-f(n)+1}}$. Then the set of all pairs of the form~$(\ell, n)$ such that~$\ell \in D_n$ is a bounded request set and we are done by applying the Kraft-Chaitin theorem.
\qed\end{pf*}

\section{Basics on Solovay functions}\label{sec:basics}


\subsection{Existence}\label{sec:existence}

\begin{theorem}[Solovay~\cite{Solovay1975}]\label{thm:existence-solf}
There exists a Solovay function.
\end{theorem}

\begin{pf*}{Proof}
Let us start by an observation. For some constant~$c$, given strings~$p$ and~$x$ such that the standard universal prefix-free machine~$\UU$ on input~$p$ outputs $x$ after exactly~$t$ steps of computation, we have $$\complK\left( \langle x,p,t \rangle \right) \leq |p|+c,$$
where $\langle .,.,. \rangle$ is the usual effective bijection from $\fs \times \fs \times \fs$ to $\fs$. This holds because there is a prefix-free machine that on input~$p$ simulates~$\UU$ on input~$p$ and outputs~$\langle x,p,t \rangle$ in case the simulated computation terminates in exactly~$t$ steps with output~$x$. Suppose now that $p$ is in addition a \emph{shortest} $\UU$-description for~$x$, i.e., $\UU(p)=x$ and $\complK(x)=|p|$. We then have
$$
|p| = \complK(x) \leq \complK\left( \langle x,p,t \rangle \right) +\bigO(1)\leq |p|+\bigO(1).
$$
Now let the function~$g_{\mathrm{S}}$ be defined by
$$
g_{\mathrm{S}}(\langle x,p,t \rangle)=
\begin{cases}
|p| &\text{if $\UU$ on input~$p$ outputs~$x$ in} \\ 
&\text{exactly $t$ steps of computations,} \\ 
2 |\langle x,p,t \rangle| &\text{otherwise}.
\end{cases}
$$
By construction, we have $\complK \leq g_{\mathrm{S}}+\bigO(1)$. Furthermore, $g_{\mathrm{S}}(\langle x,p,t \rangle)$ is at most~$\complK(\langle x,p,t \rangle)+\bigO(1)$ for all triples $\tuple{x,p,t}$ such that $p$ is a shortest $\UU$-description for~$x$ and $\UU$ on input~$p$ outputs~$x$ in exactly~$t$ steps of computation. Thus, $g_{\mathrm{S}}$ is as desired.\qed\end{pf*}

\subsection{A criterion for being a Solovay function}\label{sec:hoelzletal}

The next theorem provides the fundamental characterization of Solovay functions in terms of the sum $\sum_n 2^{-f(n)}$. It was first proven in an earlier conference paper by Bienvenu and Downey~\cite{BienvenuD2009} for computable functions, and then extended by H\"olzl, Kr\"{a}ling and Merkle~\cite{HolzlKM2013} to right-c.e.\ functions. The survey paper~\cite{BienvenuS2012} by Bienvenu and Shen gives a further extension, namely to \emph{real-valued} right-c.e.\  functions which only change finitely often at each value. 

\begin{theorem}\label{thm:fund-criterion}
Let $f \colon \N \rightarrow \N$ be a right-c.e.\ function. Then~$f$ is  a weak Solovay function if and only if $\sum_n 2^{-f(n)}$ is finite and is a Martin-L\"{o}f random real.
\end{theorem}

\label{disc-reals} There are two different ways to define the notion of a Martin-L\"of random real, however, both options lead to the same randomness notion. One option is to define Martin-L\"of randomness for infinite binary sequences first, and define a real number to be Martin-L\"of random if its binary expansion is. One other option is to directly adapt the definition of Martin-L\"of randomness to real numbers, by saying that an open set $\U$ of $\R$ is effectively open, or $\Sigma^0_1$, if it  can be written as $\bigcup_i (a_i,b_i)$, for a c.e.\ set of pairs $\{(a_i,b_i)\}_{i \in \omega}$, and considering the Lebesgue measure on $\R$ instead of the uniform measure on~$\cs$, the rest of the definition remaining the same. The second approach has been applied elegantly by Ku\v{c}era and Slaman~\cite{KuceraS2001} for proving that left-computable Martin-L\"{o}f  random reals are Solovay complete.  H\"olzl, Kr\"{a}ling and Merkle~\cite{HolzlKM2013} use a similar argument in their proof of Theorem~\ref{thm:fund-criterion}. We include this proof for ease of reference. 

\begin{pf*}{Proof of Theorem~\ref{thm:fund-criterion}}
We first show the backwards direction of the equivalence asserted in the theorem. 
We assume that $f$ is not a Solovay function and construct a sequence~$U_0, U_1, \dots$ of sets that is a Martin-L\"{o}f test and covers~$\Omega_f=\sum_n 2^{-f(n)}$.
First, we fix an appropriate computable approximation from above to~$f$, i.e., a computable function~$(n,s)\mapsto f_s(n)$ such that for all~$n$ the sequence~$f_0(n), f_1(n), ...$ is a nonascending sequence of natural numbers that converges to~$f(n)$, where we assume in addition that~$f_{s}(n) - f_{s+1}(n)$ is always equal to either~$0$ or~$1$. Let~$a_0=0$ and given~$a_i$, let~$a_{i+1}= a_i + d_i$, where~$d_i$ is obtained as follows. With some appropriate ordering of pairs understood, search for the next pair that either is of the form~$(n,0)$ or is of the form~$(n,s+1)$ where~$f_{s}(n) - f_{s+1}(n)=1$, and then let
\[
d_i=2^{-f_0(n)} \makebox[4em]{ or }d_i
= 2^{-f_{s+1}(n)}-2^{-f_s(n)} =2^{-f_s(n)},
\]
respectively. In this situation, say that the increase from~$a_i$ to~$a_{i+1}$ of size~$d_i$ \emph{occurs due to~$n$}. Furthermore, let~$b_{i}$ be the sum of all increases~$d_j$ such that~$j\le i$ and~$d_j$ and~$d_i$ occur due to the same~$n$.
By construction, all~$d_i$ and~$a_i$ are dyadic rationals and the~$a_i$ converge nondecreasingly to~$\Omega_f$.

For given~$c$, the component~$U_c$ is obtained as follows. Say that the index~$i$ is $c$-matched if~$d_i$ occurs due to~$n$ and it holds that
\begin{eqnarray}\label{dfgfghjhgertfghfgsfdsaf}
{2^{c+2}} b_{i} \leq {2^{-\complK(n)}}.
\end{eqnarray}
For every index~$i$ for which it could be verified that~$i$ is $c$-matched by approximating~$\complK$ from above, add an interval of size~$2 d_i$ to~$U_c$, where this interval either starts at~$a_i$ or at the supremum of all reals that are already in~$U_c$, whichever is larger. 

By construction,  $2^{c+2}$ times the sum of all~$d_i$ such that~$i$ is matched is at most~$\Omega_{\complK}<1$, hence the measure of~$U_c$ is at most~$2^{-(c+1)}$.
Moreover, the sets~$U_0, U_1, \ldots$ are uniformly c.e.,
hence they form a Martin-L\"{o}f test. This test covers~$\Omega_f$ because by the assumption that $f$ is not a Solovay function, and by the fact stated above that information content measures are upper bounds for $\complK$ up to an additive constant, it holds that \[\lim_{n\to\infty}(f(n)-\complK(n))=\infty.\]
Hence for every $c$ there is an index~$i(c)$ such that all~$i\ge i(c)$ will eventually become $c$-matched, where then~$U_c$ will accordingly be increased by an interval of length~$2 d_i$ to the right of~$a_{i(c)}$. The lengths of these intervals add up to~$\beta=2 (\Omega_f - a_{i(c)})$, hence~$U_c$ will cover, except for some gaps,  the open interval between~$a_{i(c)}$ and~$a_{i(c)}+\beta$, which contains~$\Omega_f$. 

It remains to show that~$\Omega_f$ cannot be a member of such a gap.
Such a gap can only occur in case an interval of size~$2 d_i$ is added to~$U_c$ in a situation where~$a_i$ is strictly larger than the supremum of all reals that are already in~$U_c$, where then the supremum of the gap is~$a_i \le \Omega_f$. Consequently, the real~$\Omega_f$ cannot be a member of the gap unless it is equal to the rational number~$a_i$ and is hence not Martin-L\"{o}f random.

Next we show the forward direction of the equivalence asserted in the theorem. 
If $f$ is a weak Solovay function, we already know by definition that $\alpha=\sum_n 2^{-f(n)}$ is finite. Let us now prove that $\alpha$ is a Martin-L\"of random real. Suppose it is not. Then for arbitrarily large~$c$ there exists~$k$ such that $\complK(\alpha \uh k) \leq k-c$ (this because of the Levin-Schnorr theorem, Theorem \ref{levin-s}). Given $\alpha \uh k$, one can effectively find some~$s$ such that
$$
\sum_{n >s} 2^{-f(n)} \leq 2^{-k}.
$$
Thus, by a standard bounded request  argument, one has $\complK(n| \alpha \uh k) \leq f(n)-k+\bigO(1)$ for all $n>s$. Thus, for all $n>s$:
$$
\complK(n) \leq f(n)+\complK(\alpha \uh k)-k+\bigO(1) \leq f(n)+(k-c)-k-\bigO(1) \leq f(n)-c-\bigO(1).
$$
And since~$c$ can be taken arbitrarily large, this shows that $\lim_{n \rightarrow +\infty} f(n)-\complK(n) = +\infty$, i.e., $f$ is not a Solovay function.
\qed\end{pf*}
The characterization of weak Solovay functions in Theorem~\ref{thm:fund-criterion} holds also when relativized to any oracle by virtually the same proof. In fact, H\"{o}lzl, Kr\"{a}ling and Merkle~\cite{HolzlKM2013} demonstrated the relativized version, as a joint generalization of, first, the already mentioned corresponding characterization of Solovay functions by Bienvenu and Downey~\cite{BienvenuD2009}, and, second, a remarkable characterization of the sequences that are weakly low for~$\complK$ by Miller~\cite{Miller2010}. For the sake of completeness, we shortly review the latter result, which now becomes a special case of Theorem~\ref{thm:fund-criterion}. Recall that the $\complK$-trivial sequences are exactly those that are low for~$\complK$, i.e., a sequence~$A$ is K-trivial if and only if $\complK^A \geq \complK - \bigO(1)$. The notion of lowness for~$\complK$ can be weakened as follows. We say that~$A$ is \emph{weakly low for~$\complK$} if $\complK^A(n) \geq \complK(n)-\bigO(1)$ for infinitely many~$n$. 

\begin{theorem}[Miller]\label{thm:miller}
The two following assertions are equivalent for all sequences~$A$.
\begin{itemize}
\item[(i)] $A$ is weakly low for~$\complK$.
\item[(ii)] $\Omega_{\complK}=\sum_n 2^{-\complK(n)}$ is Martin-L\"of random relative to~$A$.
\end{itemize}
\end{theorem}
\begin{pf*}{Proof}
The function~$\complK$ is right-c.e., hence $A$-right-c.e., and is an upper bound for~$\complK^A$ up to an additive constant. So by definition of weakly low the first statement holds if and only if~$\complK$ is a weak Solovay function relative to~$A$. The theorem is then immediate from the relativization of Theorem~\ref{thm:fund-criterion} to~$A$. \qed
\end{pf*}

Weak Solovay functions are by definition upper bounds for~$\complK$, hence always tend to infinity. We demonstrate next that indeed some Solovay functions do so nondecreasingly.

\begin{corollary}\label{cor:solovay-order}
There exists a Solovay function that is an order. 
\end{corollary}

\begin{pf*}{Proof}
Take a left-c.e.\  Martin-L\"of random real $\Omega$ and write $\Omega=\sum_n 2^{-k_n}$ where $(k_n)$ is a computable sequence of integers. Now, transform this sum as follows: rewrite each term $2^{-k_n}$ as a sum of finitely many terms $2^{-k'_n}$ where $k'_n = \max\{k_i \, |\,  i \leq n\}$. Let~$l_0, l_1, \ldots$ be the sequence of values~$k'_n$ with repetitions obtained this way and consider the function $f\colon i \mapsto l_i$. Then the function~$f$ is a computable order, and the sum~$\sum_i 2^{-f(i)}$ is equal to the Martin-L\"of random number $\Omega$, hence~$f$ is a Solovay function by Theorem~\ref{thm:fund-criterion}.  
\qed\end{pf*}

\section{Connections to  Martin-L\"of randomness}\label{sec:mlrandomness}
\subsection{A ``no-gap'' theorem for randomness}\label{sec:nogap}
In this section we will see how Solovay functions relate to Martin-L\"of randomness, via the G\'acs-Miller-Yu theorem, which is stated below as Theorem~\ref{thm:miller-yu}. To do so, we first investigate into gap phenomena for Martin-L\"of randomness, which are interesting in their own right.

By definition, a sequence is K-trivial if there is a constant~$c$ such that for all~$n$ we have that $\complK(A\uh n)\le \complK(n)+c$. Csima and Montalb{\'a}n~\cite{CsimaM2005} proved that there is an order~$h$ such that $\complK(A\uh n)\le \complK(n)+h(n)$ implies that~$A$ is $\complK$-trivial; in fact their function~$h$ is $\Delta_4^0$. Since there is a prefix-free machine that on input~$p$ outputs~$|\UU(p)|$ whenever~$\UU(p)$ is defined, there is a constant~$c$ such that for all strings~$\sigma$ it holds that~$\complK(|\sigma|) - c \le \complK(\sigma)$. Hence their result can be restated as follows: a sufficient condition for a sequence to be K-trivial is that for all~$n$ the prefix-free Kolmogorov complexity of the length~$n$ prefix of the sequence lies in the gap formed by~$K(n)-c$ and~$K(n)+h(n)$. In what follows, we will show ``no-gap theorems" for the concept of Martin-L\"{o}f randomness in the sense that for the latter concept it is not possible to obtain similarly sufficient or necessary conditions in terms of gap functions. 
  
Chaitin~\cite{Chaitin1987} proved that when a sequence~$A$ is  Martin-L\"of random, then one does not just have $\complK(A \uh n) \geq n-\bigO(1)$, but in fact $\complK(A \uh n)-n$ tends to infinity. Together with the Levin-Schnorr characterization, this shows a dichotomy: given a sequence $A \in \cs$, either $A$ is not Martin-L\"of random, in which case $\complK(A \uh n)-n$ takes arbitrarily large negative values, or~$A$ is Martin-L\"of random, in which case $\complK(A \uh n)-n$ tends to $+\infty$. This means for example that there is no sequence $A \in \cs$ such that $\complK(A \uh n) = n +\bigO(1)$. One may ask whether this dichotomy is due to a gap phenomenon, that is: is there a function~$h$ that tends to infinity, such that for every Martin-L\"of random sequence $A$, $\complK(A \uh n) \geq n+h(n)-\bigO(1)$? Similarly, is there a function~$h$ that tends to infinity such that for every sequence~$A$, $\complK(A \uh n) \geq n-h(n)-\bigO(1)$ implies that~$A$ is Martin-L\"of random? We answer both questions, as well as their plain complexity counterparts, in the negative.

\begin{theorem}\label{thm:no-gap}
There exists no function~$h\colon\N \rightarrow \N$  that tends to infinity and such that
$$
\complK(A \uh n) \geq n -h(n) - \bigO(1)
$$
is a sufficient condition for Martin-L\"of randomness of $A$. 

\noindent Similarly, there is no function~$h\colon\N \rightarrow \N$  that tends to infinity and such that
$$
\complC(A \uh n) \geq n - \complK(n)-h(n) - \bigO(1)
$$
is a sufficient condition for   Martin-L\"of randomness of $A$.
\end{theorem}
In fact, for  both statements  we  will build  counterexamples that are not even Church stochastic (see Downey and Hirschfeldt~\cite{DowneyH2010} for the definition of Church stochasticity). 
\begin{pf*}{Proof} 
First, notice that  since we want to prove this for \emph{any} function that tends to infinity, we can restrict our attention to the nondecreasing ones. Indeed, if~$h$ is a function that tends to infinity, the function $$\widetilde{h}(n)=\min\{h(i) \mid i \geq n\}$$ also tends to infinity and $\widetilde{h} \leq h$.\\

Now, assume we are in the simple case where the function~$h$ is nondecreasing and computable. A standard technique to get a nonrandom binary sequence $B$ such that $\complK(B \uh n) \geq n-h(n)-\bigO(1)$ is the following: take a Martin-L\"of random sequence $A$, and insert zeroes into $A$ in positions $h^{-1}(0),h^{-1}(1),h^{-1}(2),\ldots$. It is easy to see that the resulting sequence $B$ is not Martin-L\"of random (indeed, not even Church stochastic), and that the Kolmogorov complexity of its initial segments is as desired. This approach was refined by Merkle, Miller et al.~\cite{MerkleMNRS2006} where the authors used an insertion of zeroes on a co-c.e.\ set of positions in order to construct a left-c.e.\ sequence $B$ that is not Mises-Wald-Church stochastic, but has initial segments of very high complexity.\\

Of course, the problem here is that the function~$h$ in the hypothesis may be noncomputable, and in particular may grow slower than any computable nondecreasing function. In that case, the direct construction we just described does not necessarily work: indeed, inserting zeroes at a noncomputable set of positions may not affect  the complexity of $A$. To overcome this problem, we invoke the Ku\v{c}era-G\'acs theorem, see  Ku\v{c}era~\cite{Kucera1985}, G\'acs~\cite{Gacs1986}, or Merkle and Mihailovi\'{c}~\cite{MerkleM2004}. This theorem states that any sequence, and thus any function from $\N$ to $\N$, is Turing-reducible to some Martin-L\"of random sequence.  Hence, instead of choosing \emph{any} Martin-L\"of sequence $A$, we pick one that computes the function $h^{-1}$ and then insert zeroes into $A$ at positions $h^{-1}(0),h^{-1}(1),\ldots$ (by this we mean that we construct a sequence $B$ by adding $h^{-1}(0)$ bits of $A$, then a zero, then $h^{-1}(1)$ bits of $A$, etc.). Intuitively, the resulting sequence $B$ should not be random, as the bits of $A$ can be used to compute the places where the zeroes have been inserted. This intuition however is not quite correct, as inserting the zeroes may destroy the Turing reduction $\Phi$ from $A$ to $h^{-1}$. In other words, looking at $B$, we may not be able to distinguish the bits of $A$ from the inserted zeroes.\\

The trick to solve this last problem is to delay the insertion of the zeroes to ``give enough time'' to the reduction $\Phi$ to compute the positions of the inserted zeroes. More precisely, we insert the $k$-th zero in position $n_k=h^{-1}(k)+t(k)$ where $t(k)$ is the time needed by $\Phi$ to compute $h^{-1}(k)$ from $A$. This way, $n_k$ is computable from $A \uh n_k$ in time at most $n_k$. From this, it is not too hard to construct a computable selection rule that selects precisely the inserted zeroes, witnessing that $B$ is not Church stochastic, hence is not Martin-L\"of random. Moreover, since the ``insertion delay'' only makes the inserted zeroes more sparse, we have $\complK(B \uh n)  \geq n-h(n)-\bigO(1)$. And similarly, since $A$ is Martin-L\"of random, we have by the G\'acs-Miller-Yu theorem (Theorem~\ref{thm:miller-yu} below): $\complC(B \uh{n}) \geq n -\complK(n)-h(n)-\bigO(1)$.\\

\newcommand{\sequencex}{X} 
The formal details are as follows. Let $h$ be a nondecreasing function. By the Ku\v{c}era-G\'acs theorem, let $A$ be a Martin-L\"of random sequence and $\Phi$ be a Turing functional such that $\Phi^A(n)=h^{-1}(n)$ for all~$n$. Let~$t(n)$ be the computation time of $\Phi^A(n)$, where we can assume that~$t$ is an increasing function. Let $B \in \cs$ be the sequence obtained by inserting zeroes into~$A$ in positions $h^{-1}(n)+t(n)$. To show that $B$ is not Church stochastic, we construct a (total) computable selection rule that filters the inserted zeroes from~$B$. Let $S$ be the selection rule that works as follows on a given sequence $\sequencex \in \cs$. We proceed by induction; we call $k_n$ the number of bits selected by $S$ from $\sequencex \uh n$ and $x_n$ the prefix $\sequencex \uh n$ of $\sequencex$ from which these $k_n$ bits are deleted ($x_0$ is thus the empty string, and $k_0=0$). 

At stage~$n+1$, having already read $\sequencex \uh n$, $S$ computes $\Phi^{x_n}_n(k_n)$. If the computation halts after $s$ steps, $S$ checks whether $\Phi^{x_n}_n(k_n)+s$ returns~$n$ (the subscript of $\Phi$ refers to the number of computation steps allowed). If so, $S$ selects the $n$-th bit of $\sequencex(n)$ of $\sequencex$ and then sets $x_{n+1}=x_n$ and $k_{n+1}=k_n+1$. Otherwise, $S$ just reads the bit $\sequencex(n)$, extends~$x_n$ by this bit, i.e., $x_{n+1}=x_n\sequencex(n)$, and lets~$k_{n+1}=k_n$.\\

It is clear that $S$ is a total computable selection rule. Now suppose that we run it on~$B$. We argue that $S$ selects exactly the zeroes that have been inserted into~$A$ to get~$B$.  We prove this by induction. If $S$ has already selected from $B$ the first~$i$ inserted zeroes, then the next selected bit is the bit in position $n=\Phi^{x_n}(k_n)+s$ where $\Phi^{x_n}(k_n)$ is computed in $s$ steps. But since the selected bits are exactly the zeroes that were inserted in $A$, we have $k_n=i$ and $x_n=A \uh n-i$, and thus $s$ is the computation time of $\Phi^{x_n}(k_n)=\Phi^{A \uh n-i}(i)$, which we called~$t(i)$. And by definition of $\Phi$, $\Phi^{A \uh n-i}(i)=h^{-1}(i)$. Therefore, $n=h^{-1}(i)+t(i)$, i.e., the selected bit was an inserted zero. This proves that $S$ only selects bits that belong to the zeroes that were inserted into $A$. Conversely, we need to prove that all such bits are indeed selected by~$S$. Let $i \in \N$. The $i+1$-th inserted zero is in position $n=h^{-1}(i)+t(i)$. At stage~$n$, we have by the induction hypothesis $x_n=A \uh n-i$ and $k_n=i$. Thus, $\Phi^{x_n}_n(k_n)=\Phi^{A \uh t(i)+h^{-1}(i)-i}_{h^{-1}(i)+t(i)}(i)$, which has to halt because both quantities $t(i)+h^{-1}(i)-i$ and $h^{-1}(i)+t(i)$ are greater than $t(i)$, which is the computation time of $\Phi^A(i)$. Thus the bit in position~$n$ is indeed selected. Therefore, $S$ satisfies the desired properties, and witnesses the fact that $B$ is not Church stochastic.

Finally, for all~$n$, calling $i$ the number of inserted zeroes in $B \uh n$, we easily see that $B \uh n$ and $A \uh n-i$ can each be computed from the other one, by successive insertion or deletion of zeroes. Thus:
$\complK(B \uh n)=\complK(A \uh n-i) +\bigO(1)  \geq n-i - \bigO(1)$ since $A$ is Martin-L\"of random. And by definition of the positions where the zeroes are inserted, we have $n \geq h^{-1}(i-1)+t(i-1)$, hence $i \leq h(n)+\bigO(1)$. Therefore:
$$
\complK(B \uh n) \geq n-i-\bigO(1) \geq n-h(n) - \bigO(1)
$$
for all~$n$, which completes the proof of the first assertion. We omit the almost identical argument for the second assertion on $\complC$-complexity.
\qed\end{pf*}

Although we do not discuss them in this paper, the above construction can also be  applied to two other variants of Kolmogorov complexity, namely monotone complexity $\complKm$ and a-priori complexity $\complKM$ (see~\cite{DowneyH2010} for a definition of these complexities). A sequence~$A$ is Martin-L\"of random if and only if $\complKM(A \uh n) = \complKm(A \uh n) + \bigO(1) = n+\bigO(1)$, but there is no way to get a weaker sufficient condition. 

\begin{proposition} \label{prop:no-gap-KM}
There exists no function~$h \colon \N \rightarrow \N$ that tends to infinity and such that
$$
\complKM(A \uh n) \geq n -h(n) - \bigO(1)
$$
is a sufficient condition for $A$ to be Martin-L\"of random. Since $\complKm \geq \complKM$, this remains true with $\complKm$ in place of $\complKM$.  
\end{proposition}

The proof of Proposition~\ref{prop:no-gap-KM} is virtually identical to the proof of Theorem~\ref{thm:no-gap}. Note that the proposition is in fact stronger than the theorem, as $\complK \geq \complKM$.  \\

Another consequence of the construction performed in this proof is the dual version of Theorem~\ref{thm:no-gap} stated in Proposition~\ref{prop:no-gap-3}. The first part of the proposition has been obtained earlier on and in different ways by Miller and Yu \cite[Corollary~3.2]{MillerY2011}, and in fact with the weaker hypothesis  that $h$ is unbounded.
\begin{proposition}\label{prop:no-gap-3}
There exists no function~$h \colon \N \rightarrow \N$  that tends to infinity and such that
\[
\complK(A \uh n) \geq n + h(n) - \bigO(1)
\]
is a necessary condition for $A$ to be Martin-L\"of random. 

\noindent Similarly, there is no function~$h\colon \N \rightarrow \N$  that tends to infinity and such that
\[
\complC(A \uh n) \geq n - \complK(n) + h(n) - \bigO(1)
\]
is a necessary condition for $A$ to be Martin-L\"of random.
\end{proposition} 

\begin{pf*}{Proof}
Suppose for the sake of contradiction that there exists a function~$h'$ which tends to infinity and such that
$\complK(A \uh n) \geq n + h'(n) - \bigO(1)$ is a necessary condition for $A$ to be Martin-L\"of random. Once again, we can assume that~$h'$ is nondecreasing. Then, we perform the exact same construction as in the proof of Theorem~\ref{thm:no-gap} for a given function~$h$. At the end of proof, when evaluating the complexity of $B$, we have $\complK(B \uh n)=\complK(A \uh n-i)+\bigO(1)$, with $i \leq h(n)+\bigO(1)$, and since $A$ is Martin-L\"of random, $\complK(A \uh n-i) \geq (n-i)+h'(n-i)-\bigO(1)$. It follows that
\[
\complK(B \uh n) \geq n-h(n)+h'(n-h(n))-\bigO(1).
\]
By assumption on~$h'$, we have~$h'(n) \le n/3$ for almost all~$n$, hence if we let~$h(n)=h'(n/2)$, we have $\complK(B \uh n) \geq n-\bigO(1)$. 
This is a contradiction since by the Levin-Schnorr theorem, this would imply that the sequence~$B$ is Martin-L\"of random, which it is not by construction. The proof of the second part of the proposition is almost identical. 
\qed\end{pf*}

\subsection{The G\'acs-Miller-Yu theorem}\label{sec:gacsmilleryu}

We now turn to the link between Solovay functions and the G\'acs-Miller-Yu theorem. This theorem gives a characterization of Martin-L\"of random sequences in terms of the $\complC$-complexity of their initial segments, even though the condition still involves~$\complK$.

\begin{theorem}[G\'acs-Miller-Yu]\label{thm:miller-yu}
$A \in \cs$ is Martin-L\"of random if and only if 
\[
\complC(A \uh n) \geq n - \complK(n)-\bigO(1).
\]
Moreover there exists a computable upper bound~$f$ of~$\complK$ such that $A$ is Martin-L\"of random if and only if $\complC(A \uh n) \geq n - f(n)-\bigO(1)$.
\end{theorem}

G\'acs~\cite{Gacs1980} actually gives a variant of the first part, with conditional complexity $\complC(A \uh n|n)$ instead of $\complC(A \uh n)$. Miller and Yu~\cite{MillerY2008} proved the first part as stated above, as well as the second part about the existence of a computable~$f$ with the given properties.

The second part of the theorem indicates the existence of ``tight enough" computable upper bounds for~$\complK$. These turn out to be \emph{exactly} the Solovay functions. 

Using the ``no-gap" theorems of the previous section, we first show that any such function must be a weak Solovay function, even if we only assume the function to be merely right-c.e.

\begin{theorem}\label{thm:miller-yu-to-solovay}
Let $f $ be a right-c.e.\ function such that
\[
\complC(A \uh n) \geq n - f(n)-\bigO(1)\; \Leftrightarrow  \;
\text{$A$ is Martin-L\"of random}.
\]
Then~$f$ is a weak Solovay function. In particular, $f$ is a Solovay function in case~$f$ is computable.
\end{theorem}

\begin{pf*}{Proof}
For a start suppose that~$f$ is an upper bound for~$\complK$ up to an additive constant. Then~$f$ must be a Solovay function because otherwise $h=f-\complK$ tends to infinity and by assumption on~$f$, one has for all sequences~$A$ that
\[
\complC(A \uh n) \geq n - \complK(n) - h(n)-\bigO(1) 
\]
implies that~$A$ is Martin-L\"of random, which contradicts the no-gap result stated in Theorem~\ref{thm:no-gap}. 
In particular, in this case the function~$f$ is a Solovay function in case it is computable. (We note in passing that a similar argument shows for any, not necessarily right-c.e.\ function~$f$ that satisfies the forward implication of the equivalence in the theorem that in case the function~$f$ is an upper bound for~$K$, then this upper bound must be infinitely often tight.)   

By the preceding discussion, it suffices to show that~$f$ is an upper bound for~$\complK$ up to an additive constant. For a right-c.e.\ function, this is equivalent to~$\sum_n 2^{-f(n)} < \infty$. For the sake of contradiction, suppose that $\sum_n 2^{-f(n)} = \infty$. We already know by the G\'acs-Miller-Yu theorem that every Martin-L\"of random real~$A$ satisfies $\complC(A \uh n) \geq n - \complK(n) - \bigO(1) \geq n - 2\log (n) - \bigO(1)$.  Thus, after replacing $f(n)$ by $\min(2\log n,f(n))$ (a transformation which preserves the right-c.e.-ness and the property $\sum_n 2^{-f(n)} = \infty$), we can assume that $f(n) \leq 2\log n$. 

Let us now build an auxiliary right-c.e.\ function~$g$ as follows. Since $\sum_n 2^{-f(n)} = \infty$ and $f$ is right-c.e., one can effectively find a partition of $\N$ into consecutive intervals $I_0, I_1, I_2, \ldots$ together with stages $t_0<t_1<t_2\ldots$ such that $\sum_{n \in I_k} 2^{-f_{t_k}(n)} \geq 2^{2k}$ for all~$k$. One can further assume that $\min I_k \geq 2^k-1$ (indeed one can always add more elements to $I_{k-1}$ in order to increase $\min I_k$ if necessary). Now for all~$k$, the function~$g$ is defined on $I_k$ by $g(n)=f_{t_k}(n)+k$. This implies that $g$ is computable, that $g-f \rightarrow \infty$ and
\[
\sum_n 2^{-g(n)} = \sum_k \sum_{n \in I_k} 2^{-f_{t_k}(n)-k} \geq \sum_k 2^{k} = \infty. 
\]
Finally, observe that on every interval $I_k$, we have $g(n)=f_{t_k}(n)+k$, and by our assumption that $\min I_k \geq 2^k-1$, this means that $g(n)=f_{t_k}(n)+\bigO(\log n)=\bigO(\log n)$ (for the last inequality we use the fact that $f=\bigO(\log n)$, and assume without loss of generality that the enumeration of~$f$ from above is $\bigO(\log n)$ at all stages). This last property ensures that $\complC( n | n-g(n)) = \bigO(1)$ because the function $n \mapsto n-g(n)$ is computable and $\bigO(1)$-to-one. We can then use a well-known result due to Martin-L\"of (see~\cite[Theorem 3.11.2]{DowneyH2010}), which states that when~$g$ is a computable function such that $\complC( n | n-g(n)) = \bigO(1)$, then for any $A \in \cs$, there are infinitely many $n$ such that $\complC(A \uh n) \leq n - g(n) - \bigO(1)$. But when $A$ is Martin-L\"of random, this contradicts the hypothesis that $\complC(A \uh n) \geq n - f(n)-\bigO(1)$ (because $g-f \rightarrow \infty$). This gives us the desired contradiction, and thus proves that~$f$ is an upper bound for~$\complK$ up to an additive constant.
\qed
\end{pf*}

We now show the converse of Theorem~\ref{thm:miller-yu-to-solovay}.

\begin{theorem}\label{thm:solovay-to-miller-yu}
Let $g$ be a weak Solovay function. Then~${A \in \cs}$ is Martin-L\"of random if and only if 
\[
\complC(A \uh n) \geq n-g(n)-\bigO(1) .
\]
\end{theorem}

We begin our proof with  a combinatorial lemma.

\begin{lemma}\label{lem:memory-allocation}
Let $\sigma$ be a string. Let $I=[s,t]$ be a finite interval of integers with $s \geq |\sigma|$. Let $(a_i)_{i \in I}$ be a finite set of integers such that
\[
\sum_{i \in I} a_i 2^{-i} \geq 2^{-|\sigma|+1}.
\]
Then, there exists a subset $J$ of $I$ and a finite set of strings $S$ such that
\begin{itemize}
\item[(i)] $[S] = [\sigma]$,
\item[(ii)] for all $\tau \in S$, $|\tau| \in J$,
\item[(iii)] for all~$j \in J$, $|S \cap \{0,1\}^{\leq j}| \leq a_j$.
\end{itemize}
Moreover, $J$ and $S$ can be constructed effectively given $\sigma$, $I$ and $(a_i)_{i \in I}$.
\end{lemma}

\begin{pf*}{Proof}
We construct $J$ and $S$ via the following procedure. We initialize $J$ and $S$ to $\emptyset$. Now the procedure is as follows:\\ \\
For all $i$ from $s$ to $t$ do\\
\indent If $|S| \geq a_i $ do nothing. Otherwise:
\begin{enumerate}
\item Put $i$ into $J$.
\item Split $[\sigma] \setminus [S]$ into cylinders of measure $2^{-i}$. Let $T$ be the set of strings of length~$i$ generating those cylinders. 
\item Let $T'$ be the set containing the $c_i=a_i-|S|$ first strings of $T$ in the lexicographic order (if $c_i > |T|$ then let $T'=T$). 
\item Enumerate all strings of $T'$ into~$S$.
\end{enumerate}                                                                                                                                                                                                                                                                                                                                                                                                                                                               

We now verify that this procedure works, i.e., that the algorithm is well-defined and that the set $S$ we obtain after the $t$-loop is as wanted. First, notice that at the beginning of the $i$-loop, $S$ contains only strings of length smaller than~$i$, therefore $[S]$ can be split into cylinders of measure $2^{-i}$. Since $|\sigma| \leq s \leq i$, this is also the case for $[\sigma]$, hence for $[\sigma] \setminus [S]$, so step (2) is well-defined. We also immediately see that the conditions $(ii)$ and $(iii)$ of the lemma are satisfied: indeed, we only enumerate strings of a given length~$i$ after enumerating~$i$ into~$J$, and if we do so, we ensure that at the end of the $i$-loop, the cardinality of $S \cap \zo^{\leq i}$ is at most $a_i$. It remains to verify condition~$(i)$. First it is clear that $S \subseteq [\sigma]$ as we only enumerate cylinders that are contained in~$[\sigma]$. Suppose that this inclusion is strict. Then, when running the above procedure, at step~3, we are never in the case where $c_i > |T|$, hence for all~$i$, at the end of the $i$-loop, we have $|S \cap \zo^{\leq i}| \geq a_i$, whether~$i$ is in~$J$ or not. Therefore, at the end of the procedure, we have
\begin{align*}
\sum_{i=s}^t a_i 2^{-i} & \leq  \sum_{i=s}^t |S \cap \zo^{\leq i}| 2^{-i} 
  \leq  \sum_{i=s}^t \sum_{k=s}^i |S \cap \zo^k| 2^{-i} \\
&  \leq  \sum_{k=s}^t |S \cap \zo^k| \sum_{i=k}^t 2^{-i} 
 <   \sum_{k=s}^t |S \cap \zo^k| 2^{-k+1}\\
&  <  2 \lambda([S]) 
   <   2 \lambda([\sigma])
   <   2^{-|\sigma|+1},
\end{align*}
and this contradicts the hypothesis of the lemma. 
\qed\end{pf*}

\begin{pf*}{Proof of Theorem~\ref{thm:solovay-to-miller-yu}}
Let $g$ be a weak Solovay function. In the equivalence to be proved, 
the implication from left to right follows directly from the G\'acs-Miller-Yu theorem. In order to prove the reverse implication, let $A$ be a sequence that is not Martin-L\"of random. We shall prove that $\complC(A \uh n) \leq n-g(n)-k$ holds for infinitely many~$n$ and arbitrarily large~$k$. 

By Corollary~\ref{cor:weak-solovay-domination} below, for every  weak Solovay function~$h$ there is a Solovay function~$\widetilde{h}\le h$. The proof of the corollary does not depend on Theorem~\ref{thm:solovay-to-miller-yu} or any result demonstrated by using this theorem, hence we can apply the lemma already now and can assume that~$g$ is computable. We further assume, for technical reasons which will become clear at the end of the proof, that for all~$i$, either $g(i) \leq 2\log(i)$ or $g(i)=+\infty$. If it is not the case, replace $g$ by the bigger function $\widetilde{g}$ defined by $\widetilde{g}(i)=g(i)$ if $g(i) \leq 2\log(i)$, and $\widetilde{g}(i)=+\infty$ otherwise. Then we have
\[
\sum_i 2^{-\widetilde{g}(i)} = \sum_i 2^{-g(i)} - \sum_{\substack{i \\ g(i) \geq 2\log i}} 2^{-g(i)}, 
\]
where the third sum is a computable real number as the $i$-th term is bounded by $1/i^2$. Thus $\sum_i 2^{-\widetilde{g}(i)}$ is equal to a Martin-L\"{o}f random real minus a computable real, hence is a Martin-L\"{o}f random real and thus~$\widetilde{g}$ is still a Solovay function. 

Now, let $(\U_k)_{k \in \N}$ be a Martin-L\"of test covering~$A$ and such that $\lambda(\U_k) \leq 2^{-2k-1}$ for all~$k$. We design a procedure $(P_k)$ which for all~$k$ tries to enumerate a set of strings $S_k$ such that $[S_k]=\U_k$, with additional properties on the length of the strings it contains. We ensure that this procedure succeeds for almost all~$k$ by building an auxiliary test $\mathcal{V}_k$ which tests the randomness of $\sum_i 2^{-g(i)}$. The procedure $(P_k)$ works as follows. 

\begin{enumerate}
\item Wait for a new cylinder $[\sigma]$ to be enumerated into $\U_k$.
\item Choose a large integer~$s$, say larger than $2^N$ with $N$ larger than any integer mentioned so far in the construction (including~$k$).
\item Enumerate into $\mathcal{V}_k$ the dyadic real interval
\[
\left[ \sum_{i <s} 2^{-g(i)}, 2^{-|\sigma|+1+k}+\sum_{i <s} 2^{-g(i)}\right].
\]
\item Wait for a stage~$t$ such that
$
\sum_{i \leq t} 2^{-g(i)} > 2^{-|\sigma|+1+k}+\sum_{i <s} 2^{-g(i)}.
$
\item When this happens, we have $\sum_{i=s}^t 2^{-g(i)} > 2^{-|\sigma|+1+k}$. We then apply Lemma~\ref{lem:memory-allocation} with $a_i=2^{i-g(i)-k}$ to get a finite set of strings $S^\sigma_k$ and a finite set of integers $J_k^\sigma$ such that $[S_k^\sigma] = [\sigma]$, for all $\tau \in S_k^\sigma$, $|\tau| \in J^\sigma_k$ and for all~$j \in J_k^\sigma$, $|S^\sigma \cap \zo^{\leq j}| \leq a_j$. We then put all strings of $S^\sigma_k$ into $S_k$ and go back to step 1.
\end{enumerate}  
It is possible that for some~$k$, $(P_k)$ will at some point reach step 4 and wait there forever. We claim that this can only happen for finitely many~$k$. Indeed, for a given~$k$, we have $\lambda(\mathcal{V}_k) \leq 2^{-k}$, because whenever a cylinder $[\sigma]$ enters $\U_k$ at step 1, an interval of length $2^{-|\sigma|+1+k}$ enters~$\mathcal{V}_k$, hence $\lambda(\mathcal{V}_k) \leq 2^{k+1}\lambda(\U_k) \leq 2^{-k}$. Thus, $(\mathcal{V}_k)_{k \in \N}$ is a Martin-L\"of test. Furthermore, if the procedure for $S_k$ waits forever at some step 4, this precisely means that $\sum_{i} 2^{-g(i)}$ belongs to the dyadic interval which was put into $\mathcal{V}_k$ at step 3, and thus in that case $\sum_{i} 2^{-g(i)} \in \mathcal{V}_k$. Since~$\sum_{i} 2^{-g(i)}$ is random, it can only belong to finitely many~$\mathcal{V}_k$, hence for almost all~$k$ the procedure~$(P_k)$ never waits forever at step 4. In that case, the c.e.\ set $S_k$ it builds does satisfy $[S_k]=\U_k$ by construction.

To finish the proof, let~$k$ be such that $(P_k)$ succeeds. Since~$A$ is not Martin-L\"of random, $A$ belongs to $\U_k$, hence to $[S_k]$. This means that for some~$n$, $A \uh n$ belongs to $S_k$. To describe $A \uh n$, it suffices to describe $k$ (this can be done with $2 \log k +\bigO(1)$ bits), and its position inside $S_k$. For its position inside $S_k$, we simply describe the position of $A \uh n$ inside the $S^\sigma_k$ it belongs to, when the latter is sorted in the length-lexicographic order. By construction of $S^\sigma_k$, $n$ must be in $J^\sigma_k$ (otherwise $S^\sigma_k$ would be empty), and there are at most $a_n=2^{n-g(n)-k}$ strings of length less than or equal to~$n$ in $S^\sigma_k$, and therefore we can specify the position of $A \uh n$ inside $S^\sigma_k$ with $n-g(n)-k$ bits. Thus, our description of $A\uh n$ has total length $n-g(n)-k+2\log k+\bigO(1)$. Since~$k$ can be taken as large as wanted, this will be enough to prove the theorem, but one last thing we need to check is that this description is enough to retrieve $A \uh n$. Indeed, while we give the index of $A \uh n$ inside the $S^\sigma_k$ it belongs to, we do not describe $\sigma$ explicitly. However,~$\sigma$ can be found as follows. The description of $A \uh n$ we give has length $n-g(n)-k+2\log k+\bigO(1)$. By assumption, $g(n) \leq 2\log n$ and by construction of $S^\sigma_k$, $k \leq \log s \leq \log n$. Hence our description has length between $n-3\log n+\bigO(1)$ and $n+\bigO(1)$. Hence the length of our description gives us~$n$ with logarithmic precision. This is enough to find the string~$\sigma$ such that $A\uh n$ belongs $S^\sigma_k$ because by construction of $S_k$, if $l$ is the length of some string in $S^{\sigma'}_k$ with $\sigma' \not= \sigma$, then either $2^l < n$ or $2^n < l$, and hence either $l < n-3\log n$ or $n < l - 3 \log l$. \qed
\end{pf*}

\section{Connections to K-triviality}\label{sec:ktriviality}

\subsection{K-trivial sequences}

From their incompressibility characterization, it can be seen that the Martin-L\"of random sequences are those which have initial segments of roughly maximal Kolmogorov complexity. It is natural to ask which sequences~$A$ have initial segments of  \emph{minimal} prefix free Kolmogorov complexity $\complK(A \upharpoonright n) \le \complK(n) +\bigO(1)$. Chaitin~\cite{Chaitin1975} proved that any such sequence is computable from the halting problem, and  Solovay~\cite{Solovay1975} was able to construct such a sequence that is noncomputable and computably enumerable. The class of such sequences was further studied by Downey, Hirschfeldt, Nies and Stephan~\cite{DowneyHNS2003,Nies2005}, who called them \emph{K-trivial}. 

The $\complK$-trivial sequences  turned out to have  remarkable properties. Perhaps the most   striking fact is  that they can be characterized as the sequences that are low for Martin-L\"of randomness,  or, alternatively, as the sequences that are low for prefix-free Kolmogorov complexity. In other words, a sequence~$A$ is K-trivial if and only if Martin-L\"of randomness relativized to $A$ coincides with Martin-L\"of randomness, if and only if the prefix-free Kolmogorov complexity relativized to~$A$ is within an additive constant of  the unrelativized one. 

In this section we will show  that in the definition of the notion of K-trivial, the upper bound~$\complK(n) +\bigO(1)$ can be equivalently  replaced by any weak Solovay function, and that in fact the ability to do so characterizes the Solovay functions and the weak Solovay functions.  Using this characterization, we give                                                                                    an easy, golden-run-free proof for the fact that 
every $\complK$-trivial is Turing below a c.e.\ $\complK$-trivial. Some of the results of this section were announced in an earlier conference paper~\cite{BienvenuMN2011}.

\subsection{Solovay functions characterize K-triviality}\label{sec:solovay-to-ktriv}

In what follows, we show that weak Solovay functions can be used in place of prefix-free Kolmogorov complexity to characterize $\complK$-triviality.  This means that \complK-triviality is equivalent to $g$-triviality in the sense of  the following definition, for any weak Solovay function~$g$.
\begin{definition}\label{def:f-ktriv}
Given a function~$g\colon \N \rightarrow \N$ and an integer~$c$, a sequence~$A$ is~$g$-\emph{trivial with constant~$c$} if $\complK(A \uh n) \leq g(n) + c$ holds for all~$n$. A sequence is $g$-\emph{trivial} if it is $g$-trivial for some~$c$.  
\end{definition}
The notion of $\complK$-triviality in the sense of Definition~\ref{def:f-ktriv} coincides with the usual notion of $\complK$-triviality. However, in the usual concept of $\complK$-triviality the reference to~\complK\ is with respect to the upper bound~$\complK(n)$ but surely also to the fact that we bound the $\complK$-complexity of the initial segments of the sequence under consideration. This problem could be resolved by taking $g$-trivial as an abbreviation for $\complK$-$g$-trivial, where a sequence~$A$ is $f$-$g$-trivial if $f(A \uh n) \leq g(n) + c$ holds for some constant~$c$ and all~$n$.

We start by proving the equivalence of \complK-triviality and $g_{\mathrm{S}}$-triviality, where~$g_{\mathrm{S}}$ is Solovay's original Solovay function as constructed in the proof of Theorem~\ref{thm:existence-solf}. We will in Theorem~\ref{thm:solovay-to-ktriv} below see how to extend this equivalence to any weak Solovay function.   
\begin{theorem}\label{thm:bd-ktriv}
Let $g_{\mathrm S}$ be the Solovay function constructed by Solovay. Then a sequence is $\complK$-trivial if and only if it is $g_{\mathrm S}$-trivial.
\end{theorem}

\begin{pf*}{Proof}
One direction is easy: if $A$ is $\complK$-trivial, then $\complK(A \uh n) \leq \complK(n) + \bigO(1)$, and by definition a Solovay function is an upper bound of $\complK$ up to an additive constant, hence~$A$ is $g_{\mathrm S}$-trivial.

For the other direction, let~$A$ be $g_{\mathrm S}$-trivial for some constant~$c$. Fix~$n$. Let~$p$ be a shortest prefix description for~$n$ and let~$t$ be the running time of~$p$ on~$\UU$, i.e., $|p|=\complK(n)$ and $\UU(p)=n$ in exactly~$t$ steps. Let $m=\tuple{n,p,t}$. By definition of $g_{\mathrm S}$, we have~$g_{\mathrm S}(m)=|p|$, hence it holds that
\[
\complK(A \uh m) \leq g_{\mathrm S}(m)+c = |p| +c = \complK(n) + c .
\]
The result then follows by observing that $n$ can be retrieved from~$m$, and thus $\complK(A \uh n) \leq \complK(A \uh m) +\bigO(1) \leq \complK(n)+c+\bigO(1)$.
\qed\end{pf*}
The proof of Theorem~\ref{thm:bd-ktriv} actually shows a bit more than asserted in the theorem.  
\begin{remark}\label{rem:const-pres-1} {\rm 
The equivalence of $\complK$-triviality and $g_{\mathrm{S}}$-triviality stated in Theorem~\ref{thm:bd-ktriv} holds in the strong form that triviality constants are preserved up to an additive constant.  More precisely, there is a constant~$c_{\mathrm{S}}$ such that if a sequence is $\complK$-trivial with constant~$c$, then it is $g_{\mathrm S}$-trivial with constant~$c + c_{\mathrm{S}}$, and a similar remark holds for the reverse implication. 
}
\end{remark}
 
\begin{theorem} \label{thm:solovay-to-ktriv}
Let~$g$ be a weak Solovay function. Then any sequence is $\complK$-trivial if and only if it is $g$-trivial.
\end{theorem}

Again, the implication from $\complK$-trivial to $g$-trivial is immediate. The difficulty resides in the converse. The core of the corresponding proof is the following technical lemma, which guarantees that building a bounded request set to ensure that a sequence is $g$-trivial does not ``cost more'' (in a specific sense to be  explained below) than building a bounded request set to ensure that it is $h$-trivial for a weak Solovay function $h$. A first consequence of the lemma will be Corollary~\ref{cor:weak-solovay-domination} below, which asserts that for every weak Solovay function~$h$ there is a Solovay function~$\widetilde{h} \ge h$. The lemma and the corollary will then be applied in order to demonstrate Theorem~\ref{thm:solovay-to-ktriv}.

\begin{lemma}~\label{prop:solovay-partition} 
Let $g$ be a Solovay function, and~$h$ a weak Solovay function. There exists a positive constant~$c$ and a computable partition of~$\N$ into subsequent nonempty intervals $(I_n)_{n \in \N}$ such that for all~$n$ we have $n < \min I_n$ and
\[
2^{-g(n)} \leq 2^c \sum_{i \in I_n} 2^{-h(i)} .
\]
\end{lemma}
 
\begin{pf*}{Proof}
We will actually  use the weaker hypothesis that $g$ is computable and $\sum_n 2^{-g(n)}$ is finite.  We design a procedure which uniformly in~$p$ tries to construct  a partition $(I^p_n)_{n \in \N}$ such that $2^{-g(n)} \leq 2^p \sum_{i \in I_n} 2^{-h(i)}$. The procedure works as follows:\\ \\
For $n$ from $0$ to $\infty$ do
\begin{enumerate}
\item Let $s(p,n) \in \N$ be the least  integer $> n$ which does not belong to one of the previously constructed intervals $I^p_j$ for $j<n$. 
\item Search for  some $t > s(p,n)$ large enough to have 
\[
\sum_{i=s(p,n)}^t 2^{-h_t(i)} \geq 2^{-p} 2^{-g(n)} .
\]
\item When $t$ is found,   define $I^p_n$ to be $[s(p,n),t]$. 
\end{enumerate}

\medskip
\noindent
It is possible that  the procedure of parameter~$p$, for some $n$, never finds $t$ at Step 2. When this happens, we have by construction:
\[
\sum_{i \geq s(p,n)} 2^{-h(i)} \leq 2^{-p} 2^{-g(n)} .
\]
Hence by the Kraft-Chaitin theorem, for all~$i \geq s(p,n)$:
\begin{equation} \label{eqn:WM} 
\complK(i) \leq h(i) - p  +   \complK(p,n,s(p,n)) -   g(n) +\bigO(1) .
\end{equation} 
Using  the construction, $s(p,n)$ can be described via the pair $(p,n)$ alone, hence \begin{equation} \label{eqn:WM2} \complK(p,n,s(p,n)) \leq \complK(p,n) +\bigO(1) \leq \complK(n)+2\log p +\bigO(1).\end{equation}   Since $\sum_n 2^{-g(n)}$ is finite we have $\complK(n) \leq g(n)+\bigO(1)$.   Then \eqref{eqn:WM} and \eqref{eqn:WM2}  yield for all $i \geq s(p,n)$:
\[
\complK(i) \leq h(i) - p +2\log p + \bigO(1) .
\]
Now, recall that~$h$ is a weak Solovay function so $\complK(i) \geq h(i) +\bigO(1)$ for infinitely many~$i$. Therefore the above situation can only happen for a finite number of $p$. In other words, for all~$p$ large enough, the procedure never waits forever at step~2 and hence produces effectively a partition $(I^p_n)_{n \in \N}$ of $\N$ into intervals such that for all~$n$, and each $I^p_n = [s,t]$ we obtain as wanted
\[
2^{-p} 2^{-g(n)} \leq \sum_{i= s}^t 2^{-h_t(i)} \le \sum_{i= s}^t 2^{-h(i)} . \qed
\]
\end{pf*}

Our first application of  Lemma~\ref{prop:solovay-partition}  is the following.
\begin{corollary}\label{cor:weak-solovay-domination}
Let~$h$ be a  weak Solovay function. There exists a Solovay function~$\widetilde{h}$ such that $h \leq \widetilde{h}$. 
\end{corollary}

\begin{pf*}{Proof}
Let~$g$ be any Solovay function (for example, the one constructed in the proof of Theorem~\ref{thm:existence-solf}). By Lemma~\ref{prop:solovay-partition}, there exists a constant~$c$ and a computable partition $(I_n)_{n \in  \N}$ of $\N$ into intervals such that for all~$n$
\[
2^{-g(n)} \leq 2^c \sum_{i \in I_n} 2^{-h(i)}.
\]
Let $\widetilde{h}\colon\N \rightarrow \N$ be the function defined as follows. For a given $i$, let $I_n$ be the interval to which~$i$ belongs, and set
\[
\widetilde{h}(i)=h_t(i)~\text{ where~$t$ is the least integer s.t.\ }~2^{-g(n)} \leq 2^c \sum_{i \in I_n} 2^{-h_t(i)} .
\]
It is clear that $\widetilde{h}$ is computable and~$h\le \widetilde{h}$. Next, 
\[
\sum_i 2^{-\widetilde{h}(i)} = \sum_n \sum_{i \in I_n} 2^{-\widetilde{h}(i)}
\]
is random. Indeed, by construction for all~$n$, $   2^{-g(n)} = \bigO(\sum_{i \in I_n} 2^{-\widetilde{h}(i)})$. Hence $\sum_n 2^{-g(n)}$ is Solovay reducible to $\sum_i 2^{-\widetilde{h}(i)}$ (see \cite[p.\ 405]{DowneyH2010} for the definition of Solovay reducibility). Since the former is  random, the latter must be random as well  by the Ku\v{c}era-Slaman theorem~\cite{KuceraS2001}. Therefore~$\widetilde{h}$ is a Solovay function. 
\qed\end{pf*}
 
\begin{pf*}{Proof of Theorem~\ref{thm:solovay-to-ktriv}}
Recall that  $\UU$ is the universal prefix-free machine defining~$\complK$; a $\UU$ ``description'' of $\tau$ is  a string~$p$ such that $\UU(p)=\tau$. 

 Let~$h$ be a weak Solovay function, $d$ a constant and~$A$ a sequence such that~$\complK(A \uh n) \leq h(n)+d$ for all~$n$. We want to prove that~$A$ is $\complK$-trivial. By Corollary~\ref{cor:weak-solovay-domination}, we may assume that~$h$ is computable. We apply Lemma~\ref{prop:solovay-partition} to obtain  a constant~$c$ and a computable partition of $\N$ into intervals $(I_n)_{n \in \N}$ such that $n < \min I_n$ and  $2^{-g_{\mathrm S}(n)} \leq 2^c \sum_{i \in I_n} 2^{-h(i)}$ for all~$n$.  

We show that $A$ is $\complK$-trivial by building a bounded request set. For all~$n$ and all strings~$\sigma$ of length~$n$, we wait until we find an extension $\tau$ of $\sigma$ whose length is $\max I_n$ and such that for all~$i \in I_n$, some description of $\tau \uh i$ of length at most $h(i)+d$ is in the domain of $\UU$.  Since~$h$ is computable we can recognize when this happens. In this case, we enumerate a pair $(\sigma,g_{\mathrm S}(n)+c+d)$ into  our request set. The cost of this for us is $2^{-g_{\mathrm S}(n)-c-d}$, which we can account against the cost for  $\UU$ to enumerate descriptions of $\tau \uh i$ as above. That cost is at least $\sum_{i \in I_n} 2^{-h(i)-d}$, which in turn is at least~$2^{-g_{\mathrm S}(n)-c-d}$ by construction of the intervals $I_n$. Hence, we never spend more than~$\UU$ does, which ensures that our request set is bounded. Now, by assumption on~$A$, for every~$n$, for every $i \in I_n$, the universal machine must issue a description of $A \uh i$ of length at most $h(i)+d$, hence some pair~$(A \uh n,g_{\mathrm S}(n)+c+d)$ enters our bounded request set at some point. Therefore, for all~$n$, $\complK(A \uh n) \leq g_{\mathrm S}(n)+c+d+\bigO(1)$. By  Theorem~\ref{thm:bd-ktriv}, we can conclude  that~$A$ is $\complK$-trivial. 
\qed\end{pf*}
 
\begin{remark} \label{rem:const-pres-2}  {\rm

As in the corresponding Theorem~\ref{thm:bd-ktriv},  the equivalence of $\complK$-triviality and $g$-triviality stated in Theorem~\ref{thm:solovay-to-ktriv} holds in the strong form that triviality constants are preserved up to an additive constant.  

More precisely, every weak Solovay function~$g$ is an upper bound for~$K$ up to some additive constant~$c_g$, hence any sequence that is $\complK$-trivial with constant~$c$ is $g$-trivial with constant~$c + c_g$. Conversely, in the proof of Theorem~\ref{thm:solovay-to-ktriv} it is shown that for every weak Solovay function~$h$ there is a constant~$c$ such that $\complK(A \uh n) \leq h(n)+d$ implies $\complK(A \uh n) \leq g_{\mathrm S}(n)+c+d+\mathrm{O}(1)$, and applying Remark~\ref{rem:const-pres-1}, we get that~$A$ is $\complK$-trivial via constant~$d + c_h$ for $c_h= c+  c_{{\mathrm S}}+\mathrm{O}(1)$.} 
\end{remark} 
Another interesting corollary can be derived from the proof of Theorem~\ref{thm:solovay-to-ktriv}

\begin{remark}\label{rem:ktriv-comp-subset} {\rm
It is known that given a computable strictly increasing function $l \colon \N \rightarrow \N$, if a sequence~$X$ satisfies $\complK(A \uh l(n)) \leq \complK(n)+O(1)$, then~$X$ is \complK-trivial (see \cite[Proposition 11.1.4]{DowneyH2010}  or \cite[Exercise 5.2.9 and Solution]{Nies2009}). This fact can be extended to weak Solovay functions: For such a function~$l$, if $\complK(A \uh l(n)) \leq g(n)+O(1)$ for some weak Solovay function~$g$, then~$A$ is \complK-trivial. The proof works in the exact same way: first prove it for the particular case of the function $g$ defined on the range of~$l$ by $g(l(\langle n,p,t \rangle)) = |p|$ if $\UU(p)$ outputs~$l(n)$ in exactly~$t$ steps of computation, and outputs $2|\langle n,p,t \rangle|$ otherwise (the values outside of the range of~$l$ do not matter). The rest of the argument for the function~$g$ is the same as in Theorem~\ref{thm:bd-ktriv}). Then, extend it to all weak Solovay functions, with the same proof as Theorem~\ref{thm:solovay-to-ktriv}, only restricted to the $n$'s that are in the range of~$l$. Details are left to the reader. }
\end{remark}

\subsection{K-triviality characterizes Solovay functions}
   
Next we  prove that any right-c.e.\ function~$g$ that makes the equivalence
\begin{equation}\label{eq:ktrivialuncountable}
\text{$A$ is $\complK$-trivial} 
\makebox[7.5em]{ if and only if }
\complK(A \uh n) \leq g(n) +\bigO(1)
\end{equation}
true is a weak Solovay function, and hence a Solovay function in case~$g$ is computable. 
In the proof of our result, we need only to consider the case where~$g$ is an upper bound for~$\complK$ up to an additive constant because otherwise the class of sequences~$A$ that satisfy the right-hand side of equivalence~\eqref{eq:ktrivialuncountable} is empty. We then prove the stronger fact that in the case~$g$ is such a right-c.e.~upper bound for~$\complK$ but is not a weak Solovay function, there are uncountably many sequences~$A$ such that $\complK(A \uh n) \leq g(n) +\bigO(1)$. 
This is enough for our purposes, since there are only countably many $\complK$-trivial sequences (indeed, as we mentioned earlier, they are all computable in the halting problem). 

\begin{theorem} \label{thm:ktriv-to-solovay} 
 Suppose $g$ is  a right-c.e.\ function such that $\complK(n) \leq g(n)+\bigO(1)$ but~$g$  is not a weak Solovay   function. Then the set $\{A \mid \complK(A \uh n) \leq g(n)+\bigO(1)\}$ is uncountable.  
\end{theorem}

\begin{pf*}{Proof}
We will build an increasing sequence $a_1<a_2<a_3<\ldots$ of integers such that any subset~$A$ of $\{a_1,a_2,a_3,\ldots\}$ satisfies $\complK(A \uh n) \leq g(n)+\bigO(1)$. 

The sequence is defined  by induction (but not effectively), where we set $a_1=0$ and where we ensure by induction that for all~$k$, for any subset~$B$ of the finite set $\{a_1,\ldots,a_k\}$ and for all~$n\geq a_k$, for some constant~$d$ that depends neither on~$B$ nor on~$k$ we have that 
\begin{equation} \label{eqn:BK}  \complK(B \uh n) \le g(n) +d. \end{equation} 
This suffices to prove the desired result:  let~$A$ be any subset of $\{a_1,a_2,a_3,\ldots\}$,  and let~$n$ be  some natural number. Let~$k$ be such that $a_k\leq n<a_{k+1}$ and let $B=A \cap \{a_1,\ldots,a_{k}\}$. Then $B \uh n=A \uh n$, 
hence $\complK(A \uh n) \leq g(n) +  d$ by~\eqref{eqn:BK}.

We now explain the inductive definition of the sequence~$a_k$. Suppose we have already defined $a_1,\ldots,a_k$ with the   property (\ref{eqn:BK}). Let us choose~$c$  to be a very large integer, say~$c>2a_k+k+1$. Consider the sum $\mathrm{\Omega}_g=\sum_n 2^{-g(n)}$. By Theorem~\ref{thm:fund-criterion}, this is not a random real as $g$ is not a weak Solovay function. Hence, there exists a prefix $\sigma$ of $\mathrm{\Omega}_g$ such that $\complK(\sigma) \leq |\sigma|-c$. Let~$p$ be a shortest description for~$\sigma$. Knowing~$p$, one can effectively perform the following operations: first, retrieve~$\sigma=\UU(p)$; then, enumerate $\mathrm{\Omega}_g$ from below and wait until it becomes larger than the real value $0.\sigma$ (treated as a real number written in binary) using the approximation of the values~$g(n)$ from above; when this happens, let $a_{k+1}$ be the least number $m$~such that for all $i \geq m$, so far there has been no contribution to $\mathrm{\Omega}_g$ by the value~$g(i)$ (more precisely, via the approximation of these values from above). Since~$\sigma$ is a prefix of $\mathrm{\Omega}_g$, this means in particular that $\sum_{n\geq a_{k+1}} 2^{-g(n)}$ does not exceed $2^{-|\sigma|}$, so by the Kraft-Chaitin theorem, any integer~$n \geq a_{k+1}$ can be described by $p$ and some prefix-free code of length $g(n)-|\sigma|$. Therefore, if~$n \geq a_{k+1}$ and $B$ is a subset of $\{a_1,\ldots,a_{k+1}\}$, then $B \uh n$ can be described in a prefix-free way by
\begin{itemize}
\item[-] $B \uh {a_k}$,
\item[-] $p$ (from which $a_{k+1}$ can be retrieved),
\item[-] the single bit $B(a_{k+1})$,
\item[-] some additional $g(n)-|\sigma|$ bits. 
\end{itemize}
Thus $\complK(B \uh n) \leq 2a_k+|p|+1+g(n)-|\sigma|+\bigO(1) \leq g(n)+\bigO(1)$, using the fact that  $c>2a_k+1$ and $|p| \leq |\sigma|-c$). This concludes the inductive step. 
\qed\end{pf*}

A corollary of Theorem~\ref{thm:ktriv-to-solovay}  is that there is no $\Delta^0_2$ ``gap" for $\complK$-triviality. In the proof of this corollary we use a folklore fact about approximable orders that is stated in the following lemma.

\begin{lemma}\label{lem:delta2-order}
For every $\Delta^0_2$ order $h$ there is a right-c.e.\ order $g$ where~$g \le h$. 
\end{lemma}
\begin{pf*}{Proof}
Given an order~$h$ and a natural number~$k$, let the $k$-block of~$h$ be the finite and possible empty set of all~$n$ such that~$h(n)=k$. Observe that we have~$h_0 \le h_1$ for two given orders~$h_0$ and~$h_1$ in case for all natural numbers~$k$ the $k$-block of~$h_0$ is at least as large as the $k$-block of~$h_1$. Similarly, by increasing the size of any block of any order, the order is transformed into a strictly smaller order.

Given a $\Delta^0_2$ order $h$, write~$h$ as the pointwise limit of a uniformly computable sequence of functions $(h^0_s)$ and let
\[
h_s(n) = \max _{ i \in \{0, \ldots, n\}  } h^0_{n+s}(i) .
\]
The sequence~$(h_s)$ is uniformly computable and converges pointwise to~$h$ because 
for each~$n$ and for almost all~$s$, each of the values~$h^0_{n+s}(i)$ where~$i \le n$ agrees with~$h(i)$, and since~$h$ is an order, their maximum~$h_s(n)$ agrees with~$h(n)$. Furthermore, since the~$h^0_s$ converge to the order~$h$, for every natural number~$k$ there are~$n_0$ and~$s_0$ such that~$k <  h^0_s(n_0)$ for all~$s \ge s_0$, hence~$k < h_s(n)$ holds for all~$n \ge \max\{n_0, s_0\}$ and all~$s$. 

Let~$z_{\langle s, n\rangle}=h_s(n) $ where~$\langle .,. \rangle$ is the usual effective and effectively invertible pairing function. Then the sequence~$z_0, z_1, \ldots$ is computable and each natural number~$k$ occurs in the sequence only finitely often but at least the size of the $k$-block of~$h$ many times. Let~$g_0\colon n \mapsto n$ and for each~$i\ge 0$ let~$g_{i+1}$ be the order obtained from the order~$g_i$ by extending the size of the $z_i$-block of~$g_i$ by~$1$. Then the~$g_i$ form a uniformly computable sequence. By the preceding discussion, for all~$i$ the order~$g_{i+1}$ is strictly smaller than the order~$g_i$ and, in particular, $g_{i+1}(n) \le g_i(n)$ for all~$n$. Moreover, the~$g_i$ converge to a right-c.e.\ order~$g$. Furthermore, we have~$g \le h$ because by construction for all~$k$ the $k$-block of~$g$ is at least as large as the $k$-block of~$h$.  
\qed\end{pf*}
\begin{corollary}\label{cor:no-delta2-gap}
There is no $\Delta^0_2$ order $h$ such that for all sequences~$A$,  the sequence~$A$ is $\complK$-trivial if and only if $\complK(A \uh n) \leq \complK(n) + h(n) +\bigO(1)$.  
\end{corollary}
\begin{pf*}{Proof}
If~$h$ were right-c.e.\ this would follow directly from Theorem~\ref{thm:ktriv-to-solovay} because in this case $K+h$ would be a right-c.e.\ function that is an upper bound for~$\complK$ but is not a weak Solovay function since the order~$h$ tends to infinity. In case~$h$ is merely $\Delta^0_2$, fix a right-c.e.\ order $g \leq h$ according to Lemma~\ref{lem:delta2-order}. By the previous discussion there is a sequence~$A$ which is not $\complK$-trivial and satisfies $\complK(A \uh n) \leq \complK(n)+g(n) +\bigO(1)$ and a fortiori $\complK(A \uh n) \leq \complK(n)+h(n) +\bigO(1)$. 
\qed\end{pf*}

Corollary~\ref{cor:no-delta2-gap} was independently proven by Barmpalias and Vlek~\cite{BarmpaliasV2011}. Furthermore, Baartse and Barmpalias~\cite{BaartseB2010} showed that there \emph{does exist} a $\Delta^0_3$ gap.

\subsection{Covering $\complK$-trivials by c.e.\ $\complK$-trivials}  \label{s:deg K-trivs}

An important property of the class of $\complK$-trivial sequences is that its c.e.\ members form a Turing cover of the whole class. That is, for every $\complK$-trivial sequence~$A$, there is a c.e.\ $\complK$-trivial sequence~$B$ such that $A \leq_{\mathrm T} B$. The original proof \cite{Nies2005} is difficult and uses advanced techniques such as ``cost functions", but yields the stronger result $A \leq_{\mathrm tt} B$ (also see \cite[Corollary 5.5.3]{Nies2009}). Using Solovay functions, we obtain an elementary proof for the case of Turing reducibility, where the core of the argument relies on the following proposition.

\begin{proposition}
\label{prop_paths}
Let $A$ be a  $\complK$-trivial sequence. Then $A$ is a path of some $\complK$-trivial c.e.\ tree~$T$ which only has  finitely many paths.  
\end{proposition}
\begin{pf*}{Proof}
By Corollary~\ref{cor:solovay-order}, let~$g$ be a computable order which is also a Solovay function. For each~$k$, let $n_k=g^{-1}(k)$. Note that the sequence $(n_k)$ is computable and nondecreasing (but not necessarily increasing).  Let~$c$ be a constant such that $\complK(A \uh n) \leq g(n) +c$ for all~$n$, and $g(n) \leq \complK(n)+c$ for infinitely many~$n$. Consider the set of strings
\[
S = \{ \tau \, \mid \,  (\exists k)\,  |\tau|=n_k \, \wedge \, (\forall \sigma \leqx \tau)\, \complK(\sigma) \leq g(|\sigma|) +c \} .
\]
The set $S$ is c.e.\ and contains all the initial segments of~$A$ of type $A \uh n_k$ for some~$k$. Let~$T$ be the closure under prefixes of $S$; this makes~$T$ a c.e.\ tree such that $S$ is the restriction of~$T$ to levels of type $n_k$ for some~$k$. 

We claim that $T$ is as wanted. First of all, $A$ is a path of~$T$ by construction.  $T$ has only  finitely many paths because~$g$ is a Solovay function: any path~$B$ of~$T$ satisfies $\complK(B \uh n) \leq g(n) + c$ for all~$n$.  Hence, by Remark~\ref{rem:const-pres-2}, $B$ is $\complK$-trivial via a constant~$c+O(1)$. There are at most $2^{c+O(1)}$ such sequences~\cite{Zambella1990}. 

It remains to show that $T$ is $\complK$-trivial. By Remark~\ref{rem:ktriv-comp-subset}, we only need to prove that $\complK(T_k) \leq g(n_k) + O(1)$ for all~$k$, where $T_k$ is the restriction of $T$ to strings of length at most~$n_k$.  Fix a~$k$. $S$ being c.e.\ let $\tau$ be the last string of length at most $n_k$ enumerated into~$S$. By definition of $S$, we have $\complK(\tau) \leq g(|\tau|) +c$. Let $p$ be a description for $\tau$ of length at most $g(|\tau|)+c$ which in turn is at most $g(n_k)$ since~$g$ is nondecreasing. Up to padding~$p$ with meaningless bits, we can assume that~$p$ has length $g(n_k)+c+O(1)$. Now, given~$p$, one can retrieve $n_k$, $\tau$, the enumeration stage~$s$ of $\tau$ in $S$, all strings of~$S$ of length at most~$n_k$ (enumerating~$S$ during $s$ steps) and finally, closing under prefixes, all strings~$T$ of length at most~$n_k$. Thus
\[
\complK(T_k) \leq |p| +O(1) \leq g(n_k)+O(1) ,
\]
which is what we wanted to prove. 
\qed\end{pf*}

\begin{corollary}
Every $\complK$-trivial sequence is computable in some c.e.\ $\complK$-trivial sequence. 
\end{corollary}

\begin{pf*}{Proof}
Let $A$ be a $\complK$-trivial. By the previous proposition, let $A$ be a path of some $\complK$-trivial c.e.\ tree~$T$ with only finitely many  paths. Since a tree with finitely many paths computes all its paths, $A$ is computable from~$T$.   
\qed\end{pf*}

\section{The $c$-hitting set of a Solovay function}\label{sec:hitting}

\begin{definition}
Let $f \colon \N \rightarrow \N$ be a Solovay function and let~$c$ be an integer. The \emph{$c$-hitting set} of~$f$ is the set
\[
\HS{f,c} = \{n \, \mid \, f(n) \leq \complK(n) + c\} .
\]
\end{definition}

Note that sets of the form~$\HS{f,c}$ might be empty or finite but for a fixed Solovay function~$f$, the set~$\HS{f,c}$ is infinite for all sufficiently large~$c$.  

\begin{proposition}
Let $f$ be a Solovay function and~$c$ be an integer such that $\HS{f,c}$ is infinite. Then the set~$\HS{f,c}$ is hyperimmune and Turing-complete.
\end{proposition}

\begin{pf*}{Proof}
Suppose $\HS{f,c} = \{a_0 < a_1 < a_2 < \ldots\}$ is not hyperimmune, i.e., there is a computable function $F$ such that $a_n < F(n)$ for all~$n$. Under this assumption, we shall get a contradiction by proving that all $\complK$-trivial sequences are computable (which is not the case!). The argument is the same as Chaitin's~\cite{Chaitin1976} proof that a sequence~$A$ is computable if and only if $\complC(A \uh n) \leq \log n +\bigO(1)$. 

Let $G$ be the function defined inductively by $G(0)=0$ and $G(n+1)=F(G(n))$. Consider   the computable partition of $\N$ made of the intervals $I_n=[G(n),G(n+1)-1]$. An easy induction shows that $a_{G(n)} \in I_n$ for all~$n$, hence $\HS{f,c} \cap I_n \not=\emptyset$ for all~$n$.

Let $A$ be a noncomputable $\complK$-trivial set and, by Theorem~\ref{thm:solovay-to-ktriv}, let~$d$ be a constant such that $\complK(A \uh n) \leq f(n) +d$. Consider the c.e.\ tree
\[
T = \{\tau  \, \mid \, (\forall \sigma \leqx \tau)\, \complK(\sigma) \leq f(|\sigma|) +d\} .
\]
As in the proof of  Proposition \ref{prop_paths}, $T$ has finitely many paths, among which is~$A$. For each~$n$, let
\[
m_n = \min_{i \in I_n} \left| T \cap \zo^i \right| .
\] 
We first note that the sequence $(m_n)$ is left-c.e.\ (as $T$ is c.e.). Furthermore, it is bounded: indeed, by construction of the $I_n$, there is for each~$n$ some $i \in I_n \cap \HS{f,c}$, meaning that $f(i) \leq \complK(i)+c$. Thus any string $\tau \in  T \cap \zo^i$ is such that $\complK(\tau) \leq f(|\tau|)+d \leq \complK(|\tau|) + c + d$. By the coding theorem the number of strings of the latter type of any given length is bounded from above by some constant that depends only on~$c$ and~$d$. Let then $l = \limsup_n m_n$ and $N$ such that $m_n \leq l$ for all~$n \geq N$. Since the sequence $(m_n)$ is left-c.e.\, one can computably find two sequences $N < n_0 < n_1 < n_2 < \ldots$ and $s_0 < s_1 < s_2 < \ldots$ such that for all~$k$
\[
\min_{i \in I_{n_k}} \left| T_{s_k} \cap \zo^i \right| = l
\]
(where $T_s$ is the c.e.\ approximation of~$T$ at stage~$s$, where we assume that $T_s$ is a tree for all~$s$). 
Consider the tree
\[
T^* = \{ \tau \, \mid \, (\forall i < |\tau|)(\forall k)\,  i \in I_{n_k} \Rightarrow \tau \uh i \in T_{s_k} \} .
\]
The tree~$T^*$ is computable, and all its paths are paths of~$T$. In fact $T^*$ has exactly the same paths as~$T$, because $T^*$ coincides with $T$ on infinitely many levels (namely on each $a_k \in I_{n_k}$ such that $\left| T_{s_k} \cap \zo^i \right| = l$).  Thus $A$ is a path of $T^*$, a computable tree with finitely many paths, and therefore~$A$ is computable, a contradiction.\\ 

It remain to show that $\HS{f,c}$ is Turing complete. Observe that for any Solovay function and constant~$c$ the set $\HS{f,c}$ is co-c.e. From $\HS{f,c}$, one can compute a sequence $(n_k)_{k \in \N}$ of integers such that $\complK(n_k) \geq k$ for all~$k$. Indeed, given~$k$, since any Solovay function tends to $+\infty$, it suffices to find $n \in \HS{f,c}$ such that $f(n) \geq k+c$. Then one has $k+c \leq f(n) \leq \complK(n)+c$, thus $\complK(n) \geq k$. By a result of Kjos-Hanssen, Merkle and Stephan~\cite[Theorem 2.7]{DowneyH2010,Kjos-HanssenMS2011}, the ability to compute such a sequence is equivalent to being of diagonally noncomputable degree. By Arslanov's completeness criterion, a (co)-c.e.\ set of diagonally noncomputable degree is Turing complete. 
\qed\end{pf*}

\noindent {\bf Acknowledgement.}   Downey and Nies  were supported by the Marsden fund of New Zealand. 
We would like to thank Nan Fang for helpful discussion. Furthermore, we are grateful to the anonymous referees of \emph{Theory of Computing Systems} for their insightful remarks.

\bibliographystyle{plain}
\bibliography{SolovayFunctionsFINAL}

\end{document}